\documentclass[10pt,a4paper]{article}
\usepackage[utf8]{inputenc}
\usepackage{amsmath}
\usepackage{amsfonts}
\usepackage{amssymb}

\usepackage[utf8]{inputenc}
\usepackage{microtype}

\usepackage{amsthm}

\usepackage{verbatim}
\usepackage{resizegather}
\usepackage{tikz-cd}
\usepackage{babel} 

\usepackage[bottom]{footmisc}

\usepackage{hyperref}

\newtheorem{proposition}{Proposition}[section]
\newtheorem{corollary}{Corollary}[proposition]
\newtheorem{lemma}{Lemma}[proposition]
\newtheorem{theorem}{Theorem}
\newtheorem{remark}{Remark}[proposition]

\title{A Note on Non-Isolated Real Singularities and Links}
\author{Lars Andersen}
\date{June 2021}

\begin{document}

\maketitle
\begin{abstract}
    For analytic map germs $f: (\mathbb{R}^n, 0)\to (\mathbb{R}, 0)$ having an isolated critical value in the origin with $\dim V(f)>0$ and satisfying the transversality property of D.B. Massey we show that for $c>0$ a large enough constant, and $k\in\mathbb{N}$ a large enough natural number the local Milnor-Lê fibers of $f$ and of the isolated singularity $g=f-c(\sum_{i=1}^n x_i^2)^k$ satisfies the following. There exists a homotopy equivalence between the negative Milnor-Lê fiber of $f$, to which a cobordism between its boundary and the link of $f$ has been adjoined, and the Milnor-Lê fiber of $g$. 
\end{abstract}
\section{Introduction}
\subsection{Introduction}
The topological nature of analytic map germs has grown to this day to becomee an immense field of research. For real analytic maps, and more particularly singular such maps, we can for the purpose of this paper mention the works of H. Hamm (e.g \cite{hamm}), Z. Szafraniec (\cite{Sza}) and more recently the paper of Séade, Cisneros-Molina and Snoussi (\cite{CSS}) and the preprint (\cite{dutertre2019topology}) of N. Dutertre. The result presented in following paper is in fact a mere consequence of the work of the latter authors, and should ideally be called a corollary of \cite[Theorem 13.5]{JS} and \cite[Lemma 1]{Sza}.\\\\

The result in this paper is a topological relation between the link of a singularity of real analytic function germs $f: (\mathbb{R}^n, 0)\to (\mathbb{R}, 0)$, its Milnor-Lê fiber, and the Milnor fiber of an isolated singularity of the form $f-c\omega^k$ where $\omega$ is the sum of the squares of the coordinate functions on $\mathbb{R}^n$, $c>0$ a constant, and $k\in\mathbb{N}$ a sufficiently large natural number.

\subsection{The Local Milnor-Lê Fibration}
Suppose that the real analytic function germ $f: (\mathbb{R}^n, 0)\to (\mathbb{R}, 0)$ has an isolated critical value in the origin, satisfies the transversality property of D.B. Massey (\cite[Definition 13.4]{Seade}) and that $\dim V(f)>0.$ Then the following theorem holds
\begin{theorem}[{\cite[Theorem 13.5]{JS}}]\label{one} The germ $f$ has a local Milnor-Lê fibration
$$f: N_f(\epsilon,\delta)\to \mathbb{R}\setminus\{0\},$$
where $N_f(\epsilon, \delta)=f^{-1}((-\epsilon, \epsilon)\setminus\{0\})\cap \mathbb{B}_{\delta}$ for some ball $\mathbb{B}_{\delta}\subset \mathbb{R}^n$ and for $0<\epsilon<<\delta$. This determines an equivalent fiber bundle
$$\phi: \mathbb{S}_{\delta}\setminus K\to \mathbb{S}^1$$
where $K=f^{-1}(0)\cap \mathbb{S}_{\delta}$ is the \emph{link} of the map germ $f$, and where the projection map $\phi$ is $f/\lVert f\rVert$ when restricted to $N_f(\epsilon, \delta)$.
\end{theorem}

Let $p>0$ be a positive integer. Then the transversality property is in fact equivalent to the existence of the first fibration for locally surjective map germs 
$$g: (\mathbb{R}^n, 0)\to (\mathbb{R}^p, 0),\qquad \dim V(g)>0$$
with an isolated critical value in the origin.\\

In our case ($p=1$) the transversality property implies the equivalence of the bundles. This is done, following Milnor, by integrating an appropriate vector field to "inflate" the Milnor tube $N_f(\epsilon,\delta)$ and push the fibers of the first fibration outwards towards the sphere. For $p>1$, the equivalence of the bundles is more delicate yet has recently been shown (\cite{CSS}) by J. L. Cicneros-Molina, J. Séade and J. Snoussi that the notion of $d$-regularity which they introduced, plays an essential rôle.\\\\
\subsection{A Result of Szafraniec}
Let us now recall a result of Szafraniec.

\begin{lemma}[{\cite[Lemma 1]{Sza}}]\label{Szalemma} Let $f: (U, 0)\to (\mathbb{R}, 0)$ be a  real analytic function defined in an open subset of $\mathbb{R}^n$. Then there exist constants $C>0,\alpha>0$ such that:\\ if $c\in (0, C), k\geq \alpha$ is an integer, $r\neq 0$ is sufficiently close to the origin and 
$$g: U\subset\mathbb{R}^n\to \mathbb{R},\qquad g=f-c(x_1^2+\dots+ x_n^2)^k$$
then $\{0\}\subset\mathbb{R}$ is a regular value of $g_{|\mathbb{S}_r}$. In particular $g$ has an isolated critical point in the origin.
\end{lemma}

In particular it follows from the proof of the lemma that if $r=r(c, k)$ is chosen sufficiently small and if
$$N_f^{-}(r):=\{x\in \mathbb{S}_r: f(x)\leq 0\}$$
$$N_g^{-}(r):=\{x\in\mathbb{S}_r: g(x)\leq 0\}$$
then $N_f(r)\subset \text{int} N_g(r)$ and $N_f(r)$ is a deformation retract of $N_g(r)$.\\

\subsection{The Result}
Let us from now on assume that $f: (\mathbb{R}^n, 0)\to (\mathbb{R}, 0)$ satisfies the hypotheses of Theorem \hyperref[]{\ref*{one} }. We then apply the theorem to $-f$ and consequently let $N_f^{-}(\epsilon, \delta)$ be the negative Milnor tube of $f$, so that
$$f: N_f^{-}(\epsilon, \delta)\to (-\infty, 0)$$
is the projection of a trivial fibre bundle. Let 
$$\mathcal{F}^{-}(f)=f^{-1}(\epsilon)\cap \mathbb{B}_{\delta}.$$
By inflating the Milnor tube as in the proof of Theorem \hyperref[]{\ref*{one} } one obtains a homeomorphism
$$\mathcal{F}^{-}(f)\cong \{x\in \mathbb{S}_{\delta}:\quad f(x)\leq -\epsilon\}.$$
As a consequence, 
$$\mathcal{F}^{-}(f)\cup_{\partial \mathcal{F}^{-}(f)} (\overline{N_f^{-}(\epsilon, \delta)}\cap \mathbb{S}_{\delta})\cong N_f(\delta)$$
where
$$\overline{N_f^{-}(\epsilon, \delta)}\cap \mathbb{S}_{\delta}=f^{-1}([-\epsilon, 0])\cap \mathbb{S}_{\delta}.$$
On the other hand, by Lemma \hyperref[Szalemma]{\ref*{Szalemma} }, $N_f(\delta)\hookrightarrow N_g(\delta)$ is a deformation retract so 
$$\mathcal{F}^{-}(f)\cup_{\partial \mathcal{F}^{-}(f)} (\overline{N_f^{-}(\epsilon, \delta)}\cap \mathbb{S}_{\delta})\hookrightarrow N_g^{-}(\delta)$$
is a homotopy equivalence. Note that by Lemma \hyperref[Szalemma]{\ref*{Szalemma} }, $g$ has an isolated critical point in the origin.\\

 Let us now recall a result due to A. Durfee. Suppose given an algebraic set $M$ in the affine space $\mathbb{R}^n$ and let $X\subset M$ be a compact algebraic subset. Suppose that $M\setminus X$ is nonsingular. Recall that a subset $T\subset M$ is an \emph{algebraic neighborhood} of $X$ if:\\
there exists a nonnegative proper polynomial map $\alpha: M\to \mathbb{R}$ and a positive real number $\gamma$ smaller than any critical value of $\alpha$ such that
$$X=\alpha^{-1}(0),\qquad T=\alpha^{-1}([0, \gamma])$$ 
In this situation Durfee (\cite{durfee}) proved 
\begin{lemma}[{\cite[Proposition 1.6]{durfee}}]\label{durfeelemma} If $T$ is an algebraic neighborhood of $X$ then the inclusion $X\hookrightarrow T$ is a homotopy equivalence.

\end{lemma}
\begin{remark} The main point of the proof is to pick a well-chosen strictly increasing neighborhood basis of $X$ and then use the vector field $\text{grad} \alpha$ to trivialise these.
\end{remark}

We now apply this with 
$$X=K_f=g^{-1}(0)\cap \mathbb{S}_{\delta},\quad T=g^{-1}([-\epsilon', \epsilon'])\cap \mathbb{S}_{\delta}$$
for $\epsilon'<<\delta$, and assume that 
$$g^{-1}([-\epsilon', 0))\cap \mathbb{S}_{\delta}\neq\emptyset.$$
Of course this is always the case whenever the boundary of the negative Milnor fiber $\mathcal{F}^{-}(g)$ of $g$ is nonempty. Then \hyperref[Szalemma]{Lemma \ref*{Szalemma} } together with \hyperref[one]{Theorem \ref*{one} } applied to the isolated singularity $g$ give 
$$ N_g^{-}(\delta)\sim \{g\leq -\epsilon'\}\cap \mathbb{S}_{\delta}\cong \mathcal{F}_{\eta'}^{-}(g).$$
\newpage
We have proven

\begin{theorem}\label{main} Suppose that $f: (\mathbb{R}^n, 0)\to (\mathbb{R}, 0)$ has an isolated critical value in the origin, satisfies the transversality property and that $\dim V(f)>0.$ Let $\mathcal{F}_{\epsilon}^{-}(f)$  denote the local negative Milnor-Lê fiber of $f$, where $0<\epsilon<<\delta$. There exist constants $C>0, \alpha>0$ such that:\\
if $c\in (0, C)$, if $k\geq \alpha$ is an integer and if $\delta$ is chosen so small that 
$$g: \mathbb{B}_{\delta}\to \mathbb{R},\qquad g=f-c(x_1^2+\dots+x_n^2)^k$$
has no critical points except an isolated critical point in the origin then there exists a homotopy equivalence 
$$\mathcal{F}_{\epsilon}^{-}(f)\cup_{\partial \mathcal{F}_{\epsilon}^{-}(f)} (\overline{N_f^{-}(\epsilon, \delta)}\cap \mathbb{S}_{\delta})\sim \mathcal{F}_{\epsilon'}^{-}(g)$$
whenever the negative Milnor-Lê fiber $\mathcal{F}_{\epsilon'}^{-}(g)$ of $g$ at the origin is non-empty.
\end{theorem}

\begin{corollary} Under the same assumptions as in \hyperref[main]{Theorem \ref*{main} } there is a long exact sequence in homology
$$\to H_n(\partial \mathcal{F}_{\epsilon}^{-}(f))\to H_n(\mathcal{F}_{\epsilon}^{-}(f))\oplus H_n(N^{-}(\partial \mathcal{F}_{\epsilon}^{-}(f)))\to H_n(\mathcal{F}_{\epsilon'}^{-}(g))\to $$
\end{corollary}

\bibliographystyle{plain}
\bibliography{biblio}


\end{document}